\documentclass[12pt]{article}
\usepackage{amssymb,amsmath,latexsym, verbatim}
\allowdisplaybreaks

\begin{document}

\begin{doublespace}

\newtheorem{thm}{Theorem}[section]
\newtheorem{lemma}[thm]{Lemma}
\newtheorem{defn}[thm]{Definition}
\newtheorem{prop}[thm]{Proposition}
\newtheorem{corollary}[thm]{Corollary}
\newtheorem{remark}[thm]{Remark}
\newtheorem{example}[thm]{Example}
\numberwithin{equation}{section}
\def\ee{\varepsilon}
\def\qed{{\hfill $\Box$ \bigskip}}
\def\NN{{\cal N}}
\def\AA{{\cal A}}
\def\MM{{\cal M}}
\def\BB{{\cal B}}
\def\CC{{\cal C}}
\def\LL{{\cal L}}
\def\DD{{\cal D}}
\def\FF{{\cal F}}
\def\EE{{\cal E}}
\def\QQ{{\cal Q}}
\def\RR{{\mathbb R}}
\def\R{{\mathbb R}}
\def\L{{\bf L}}
\def\K{{\bf K}}
\def\S{{\bf S}}
\def\A{{\bf A}}
\def\E{{\mathbb E}}
\def\F{{\bf F}}
\def\P{{\mathbb P}}
\def\N{{\mathbb N}}
\def\eps{\varepsilon}
\def\wh{\widehat}
\def\wt{\widetilde}
\def\pf{\noindent{\bf Proof.} }
\def\beq{\begin{equation}}
\def\eeq{\end{equation}}
\def\dint{\displaystyle\int}

\title{\Large \bf Trace estimates for relativistic stable processes}
\author{{\bf Hyunchul Park} \quad and \quad {\bf Renming Song}\thanks{Research supported in part by a grant from the Simons
Foundation (208236).}
}

\date{}
\maketitle

\begin{abstract}
In this paper, we study the asymptotic behavior, as the time $t$ goes to zero,
of the trace of the semigroup of a killed relativistic $\alpha$-stable process
in bounded $C^{1,1}$ open sets and bounded Lipschitz open sets.
More precisely, we establish the asymptotic expansion in terms of $t$ of the trace
with an error bound of order $t^{2/\alpha}t^{-d/\alpha}$ for $C^{1,1}$ open sets and
of order $t^{1/\alpha}t^{-d/\alpha}$ for Lipschitz open sets. Compared with the corresponding
expansions for stable processes, there are more terms between the orders
$t^{-d/\alpha}$ and $t^{(2-d)/\alpha}$ for $C^{1,1}$ open
sets, and, when $\alpha\in (0, 1]$, between the orders $t^{-d/\alpha}$ and $t^{(1-d)/\alpha}$ for
Lipschitz open sets.
\end{abstract}

\section{Introduction and statement of the main results}\label{Introduction and statement of the main result}
For any $m>0$ and $\alpha\in (0, 2)$, a relativistic $\alpha$-stable process
$X^{m}$ on $\R^{d}$ with mass $m$ is a L\'evy process with characteristic function given by
\beq\label{char:RSP}
\E\left[\exp(i \xi\cdot(X_{t}^{m}-X_{0}^{m}))\right]=\exp(-t((|\xi|^{2}+m^{2/\alpha})^{\alpha/2}-m)), \quad
\xi \in \R^{d}.
\eeq
The limiting case $X^{0}$, corresponding  to $m=0$, is a (rotationally) symmetric $\alpha$-stable process on $\R^{d}$ which we will simply denote as $X$.
The infinitesimal generator of $X^{m}$ is $m-(m^{2/\alpha}-\Delta)^{\alpha/2}$.
Note that when $m=1$, this infinitesimal generator reduces to $1-(1-\Delta)^{\alpha/2}$. Thus the 1-resolvent kernel of the relativistic $\alpha$-stable
process $X^{1}$ on $\R^{d}$ is just the Bessel potential kernel. When $\alpha=1$, the infinitesimal generator reduces to
the so-called free relativistic Hamiltonian $m-\sqrt{-\Delta+m^{2}}$. The operator $m-\sqrt{-\Delta+m^{2}}$ is
very important in mathematical physics due to its application to relativistic quantum mechanics.

In this paper, we will be interested in the asymptotic behavior of the trace of the semigroup
associated with killed relativistic $\alpha$-stable processes in open sets of $\R^{d}$. The process $X^{m}$ has a transition density $p^{m}(t,x,y)=p^{m}(t,y-x)$ given by the inverse Fourier transform
$$
p^{m}(t,x)=(2\pi)^{-d}\int_{\R^{d}}e^{-i\xi x}e^{-t(|\xi|^{2}+m^{2/\alpha})^{\alpha/2}+mt}d\xi.
$$
For any open set $D$ in $\R^{d}$, the killed relativistic $\alpha$-stable process $X_{t}^{m,D}$ is defined by
$$
X_{t}^{m,D}=
\begin{cases}
X_{t}^{m}  &\mbox{if } t < \tau_{D}^{m},\\
\partial &\mbox{if } t \geq \tau_{D}^{m},\\
\end{cases}
$$
where $\tau_{D}^{m}=\inf\{t>0 : X_{t}^{m} \notin D\}$ is the first exit time of $X^{m}$ from $D$.
The process $X_{t}^{m,D}$ is a strong Markov process with a transition density $p_{D}^{m}(t,x,y)$ given by
$$
p_{D}^{m}(t,x,y)=p^{m}(t,x,y)-r_{D}^{m}(t,x,y),
$$
with
$$
r_{D}^{m}(t,x,y)=\E_{x}\left[t>\tau_{D}^{m} ; p^{m}(t-\tau_{D}^{m}, X_{\tau_{D}^{m}}^{m},y)\right].
$$

We denote by $(P_{t}^{m,D}: t\ge 0)$ the semigroup of $X_{t}^{m}$ on $L^{2}(D)$: for any $f\in L^{2}(D)$,
$$
P_{t}^{m,D}f(x):=\E_{x}\left[f(X_{t}^{m,D})\right]=\int_{D}f(y)p_{D}^{m}(t,x,y)dy.
$$
Whenever $D$ is of finite volume, $P_{t}^{m,D}$ is a Hilbert-Schmidt operator
mapping $L^{2}(D)$ into $L^{\infty}(D)$ for every $t>0$.
By general operator theory, there exist
an orthonormal basis of eigenfunctions $\{\phi^{(m)}_{n}\}_{n=1}^{\infty}$ for $L^{2}(D)$ and corresponding
eigenvalues $\{\lambda^{(m)}_{n}\}_{n=1}^{\infty}$
of the generator of the semigroup $P_{D}^{m,D}$ satisfying
$$
0<\lambda^{(m)}_{1}<\lambda^{(m)}_{2}\leq \lambda^{(m)}_{3}\leq \cdots
$$
with $\lambda^{(m)}_{n}\rightarrow \infty$.
By definition, we have
$$
P_{t}^{m,D}\phi^{(m)}_{n}(x)=e^{-\lambda^{(m)}_{n}t}\phi^{(m)}_{n}(x),\quad x\in D, t>0.
$$
We also have
$$
p_{D}^{m}(t,x,y)=\sum_{n=1}^{\infty}e^{-\lambda^{(m)}_{n}t}\phi^{(m)}_{n}(x)\phi^{(m)}_{n}(y).
$$
$\lambda^{(0)}_n$ will be simply denoted by $\lambda_n$.

In the remainder of this paper, we assume $d\ge 2$.
We are interested in finding the asymptotic behavior, as $t\to 0$, of the trace defined by
$$
Z_{D}^{m}(t)=\int_{D}p^{m}_{D}(t,x,x)dx=\sum_{n=1}^{\infty}e^{-\lambda^{(m)}_{n}t}\int_{D}(\phi^{(m)}_{n})^{2}(x)dx
=\sum_{n=1}^{\infty} e^{-\lambda^{(m)}_{n}t}.
$$
It is shown in \cite{BK} that for any open set $D$ of finite volume, it holds that
\begin{equation}\label{e:new}
\lim_{t\rightarrow 0}t^{d/\alpha}Z_{D}^{0}=C_{1}|D|, \qquad C_{1}=\frac{\omega_{D}\Gamma(d/\alpha)}{(2\pi)^{d}\alpha}.
\end{equation}
This is closely related to the growth of the eigenvalues of $P_{t}^{0,D}$:
if $N^{0}(\lambda)$ is the number of eigenvalues $\lambda_j$ such that $\lambda_j\le \lambda$,
then it follows from the
classical Karamata Tauberian theorem (see for example \cite{F}) that
\beq\label{1st term : SSP}
N^{0}(\lambda) \sim \frac{C_{1}|D|}{\Gamma(d/\alpha +1)}\lambda^{d/\alpha}, \quad \text{as } \lambda \rightarrow \infty.
\eeq
This is the analogue for killed stable processes of the celebrated Weyl's 
asymptotic formula for the eigenvalues of the Dirichlet Laplacian.
We will see later in this paper that exactly the same formula is true 
for relativistic stable processes. That is, the first term in the expansion
of $Z_{D}^{m}(t)$ is the same as that of $Z_{D}^{0}(t)$ and \eqref{1st term : SSP} is also true for relativistic stable processes.

Our main goal in this paper is to get the asymptotic expansion of $Z_{D}^{m}(t)$ as $t\rightarrow 0$
under some additional assumptions on the smoothness of the boundary of $D$.
Our work is inspired by the paper \cite{B} for Brownian motion and the papers \cite{BK, BKS} for stable processes.
The first theorem is an asymptotic expansion of $Z_{D}^{m}(t)$
with error bound of order $t^{2/\alpha}t^{-d/\alpha}$ in $C^{1,1}$ open sets.
To state the result precisely, we need some definitions.
Recall that an open set $D$ in $\R^{d}$ is said to be a (uniform) $C^{1,1}$ open set
if there are (localization radius) $R>0$
and $\Lambda_{0}$ such that for every $z\in\partial D$, there exist
a $C^{1,1}$ function $\phi=\phi_{z}:\R^{d}\rightarrow \R$ satisfying $\phi(0, \cdots, 0)=0$,
$\nabla\phi(0)=(0,\dots,0)$,
$|\nabla\phi(x)-\nabla\phi(y)|\leq \Lambda_{0}|x-z|$ and an orthonormal coordinate system $CS_{z}:y=(y_{1},\cdots,y_{d-1},y_{d}):=(\tilde{y},y_{d})$ with origin at $z$
such that $B(z,R)\cap D =\{y= (\tilde{y},y_{d})\in B(0,R) \text{ in } CS_{z} : y_{d}>\phi(\tilde{y})\}$.
For $x\in\R^{d}$, let $\delta_{D}(x)$ denote the Euclidean distance between $x$ and $D^{c}$ and $\delta_{\partial D}(x)$ the Euclidean distance between $x$
and $\partial D$. It is well known that a $C^{1,1}$ open set $D$ satisfies both the $\textit{uniform interior ball condition}$ and the $\textit{uniform exterior ball condition}$:
there exists $r_{0}<R$ such that for every $x\in D$ with $\delta_{\partial D}(x)<r_{0}$ and $y\in \R^{d}\setminus \bar{D}$ with $\delta_{\partial D}(y)<r_{0}$, there are
$z_{x},z_{y}\in\partial D$ so that $|x-z_{x}|=\delta_{\partial D}(x)$, $|y-z_{y}|=\delta_{\partial D}(y)$ and that $B(x_{0},r_{0})\subset D$
and $B(y_{0},r_{0})\subset \R^{d}\setminus\bar{D}$, where $x_{0}=z_{x}+r_{0}(x-z_{x})/|x-z_{x}|$ and $y_{0}=z_{y}+r_{0}(y-z_{y})/|y-z_{y}|$.
In fact, $D$ is a $C^{1,1}$ open set if and only if $D$ satisfies the uniform interior ball condition and the uniform exterior ball condition (see \cite[Lemma 2.2]{AKSZ}).
In this paper we call the pair $(r_{0},\Lambda_{0})$ the characteristics of the $C^{1,1}$ open set $D$.
For any open set $D$ in $\R^d$, we use $|D|$ to denote the $d$-dimensional Lebesgue measure
of $D$ and $\mathcal{H}^{d-1}(\partial D)$ to denote the $(d-1)$-dimensional Hausdorff measure
of $\partial D$. When $D$ is a $C^{1, 1}$ open set, $\mathcal{H}^{d-1}(\partial D)$ is equal to
the surface measure $|\partial D|$ of $\partial D$.
We will use $H$ to denote the half space $\{x=(x_{1},x_{2},\cdots,x_{d}) : x_{1}>0\}$.

The following is the the first main result of this paper.

\begin{thm}\label{RSP:main}
Suppose that $D$ is a bounded $C^{1,1}$ open set in $\R^{d}$.
Let $k$ be the largest integer such that $k< \frac{2}{\alpha}$.
Then the trace $Z_{D}^{m}(t)$ admits the following expansion
$$
Z_{D}^{m}(t)=C_{1}|D|t^{-\frac{d}{\alpha}}-C_{2}|\partial D|t^{\frac{1-d}{\alpha}}+
\frac{\omega_{d}\Gamma(d/\alpha)|D|}{(2\pi)^{d}\alpha}t^{-\frac{d}{\alpha}}
\sum_{n=1}^{k}\frac{m^{n}}{n!}t^{n}+O(\frac{t^{2/\alpha}}{t^{d/\alpha}}),
$$
where $C_1$ is given in \eqref{e:new} and
$$
C_{2}=\int_{0}^{\infty}r_{H}^{0}(1,(r,\tilde{0}),(r,\tilde{0}))dr.
$$
\end{thm}

The second main result of the paper is an asymptotic expansion of $Z_{D}^{m}(t)$ with error bound of
order $t^{1/\alpha}t^{-d/\alpha}$ in Lipschitz open sets.
Before we state the second main result, we recall the definition
of Lipschitz open sets. An open set $D$ in $\R^d$ is called a Lipschitz open set
if there exist constants $R_{0}$ (localization radius)
and $\lambda>0$ (Lipschitz constant) such that for every $z\in\partial D$
there exist a Lipschitz function $F:\R^{d-1}\rightarrow \R$ with Lipschitz constant $\lambda$
and an orthornormal coordinate
system $y=(y_{1},\cdots,y_{d})$ such that
$
D\cap B(z,R_{0})=\{y:y_{d}>F(y_{1},\cdots,y_{d-1})\}\cap B(z,R_{0}).
$
Here is the second main result.
\begin{thm}\label{RSP:main2}
Suppose that $D$ is a bounded Lipschitz open set in $\R^{d}$.
Let $j$ be the largest integer such that $j\leq \frac{1}{\alpha}$.
Then the trace $Z_{D}^{m}(t)$ admits the following expansion
$$
t^{d/\alpha}Z_{D}^{m}(t)=C_{1}|D|-C_{2}\mathcal{H}^{d-1}(\partial D)t^{1/\alpha}+
\frac{\omega_{d}\Gamma(d/\alpha)|D|}{(2\pi)^{d}\alpha}\sum_{n=1}^{j}\frac{m^{n}}{n!}t^{n}+o(t^{1/\alpha}),
$$
where $C_1$ and $C_2$ are the same as in Theorem \ref{RSP:main}.
\end{thm}

The asymptotic behaviors of the trace $Z_D(t)$ of the killed Brownian motion
(i.\/e., killed symmetric $\alpha$-stable
process with $\alpha=2$) in bounded domains $D$ of $\R^{d}$
have been extensively studied by many authors.
It is shown in \cite{vB} that, when $D$ is a bounded $C^{1,1}$ domain,
$$
\left|Z_{D}(t)-(4\pi t)^{-d/2}\left(|D|-\frac{\sqrt{\pi t}}{2}|\partial D|\right)\right|\leq \frac{c|D|t^{1-d/2}}{R^{2}}, \quad t>0.
$$
The following asymptotic result
\beq\label{Trace BM:C1}
Z_{D}(t)=(4\pi t)^{-d/2}\left(|D|-\frac{\sqrt{\pi t}}{2}|\partial D|+o(t^{1/2})\right), \quad t\rightarrow 0,
\eeq
was proved in \cite{BC} when $D$ is a bounded $C^{1}$ domain.
\eqref{Trace BM:C1} was subsequently extended to Lipschitz domains in \cite{B}.

The asymptotic behaviors of the trace $Z_{D}^{0}(t)$ of killed
symmetric $\alpha$-stable processes, $0<\alpha<2$, in open sets of
$\R^{d}$ have been  studied in \cite{BK, BKS}.
It was shown in \cite{BK} that, for any bounded $C^{1,1}$ open set $D$,
$$
\left|Z_{D}^{0}(t)-\frac{C_{1}|D|}{t^{d/\alpha}}+\frac{C_{2}|\partial D|t^{1/\alpha}}{t^{d/\alpha}}\right|
\leq \frac{c|D|t^{2/\alpha}}{r_{0}^{2}t^{d/\alpha}},
$$
where $C_1$ and $C_2$ are the same as in Theorem \ref{RSP:main} and $c$ is a positive constant depending on $d$ and $\alpha$ only.
It was shown in \cite{BKS} that, when $D$ is a  bounded Lipschitz domain, $Z_{D}^{0}(t)$ satisfies
$$
t^{d/\alpha}Z_{D}^{0}(t)=C_{1}|D|-C_{2}\mathcal{H}^{d-1}(\partial D)t^{1/\alpha}+o(t^{1/\alpha}).
$$

\begin{remark}
Note that the first term in the expansion of $Z_{D}^{m}(t)$ is exactly the same as in the case of $Z_{D}^{0}(t)$. However the rest of the terms are quite different.
We note here that the coefficient of the term of order $t^{1/\alpha}t^{-d/\alpha}$ is the same
in the stable process case, but in the case of relativistic stable processes for $C^{1,1,}$ open sets, there are $k$ intermediate terms
of the form $t^{k}t^{-d/\alpha}$,
where $k$ is a positive integer less than $2/\alpha$.
Since $0<\alpha<2$,
there is at least one more term involved in the asymptotic expansion
of $Z_{D}^{m}(t)$ than that of $Z_{D}^{0}(t)$ up to order of $t^{2/\alpha}t^{-d/\alpha}$.
For Lipschitz open sets, when $\alpha\leq 1$ there are $j$ intermediate terms
of the form $t^{j}t^{-d/\alpha}$, where $j$ is an integer that is less than or equal to $1/\alpha$.
\end{remark}

\begin{remark}
In \cite{BMN}, an asymptotic expansion for the trace of relativisitic $\alpha$-stable processes in
bounded $C^{1, 1}$ open sets was established. Compared with Theorem \ref{RSP:main},
the expansion of  \cite{BMN}
does not contain the intermediate terms.
\end{remark}

The rest of the paper is organized as follows. In Section \ref{Preliminaries}, we recall some basic facts about relativistic stable processes and
present several preliminary results which will be used in Sections \ref{Proof of the main result}
and \ref{main result 2:Lipschitz}.
Theorem \ref{RSP:main} is proved in Section \ref{Proof of the main result},
while Theorem \ref{RSP:main2} is proved in Section \ref{main result 2:Lipschitz}.

Throughout this paper, we will use $c$ to denote a positive constant depending
(unless otherwise explicitly stated) only on $d$ and $\alpha$ but whose
value may change from line to line, even within a single line. In this paper,
the big O notation $f(t)=O(g(t))$ always means that there exist constants
$C$ and $t_{0}>0$ such that $f(t)\leq C g(t)$ for all $0< t < t_{0}$.


\section{Preliminaries}\label{Preliminaries}

In this section, we recall some basic facts about relativistic $\alpha$-stable processes. From
\eqref{char:RSP}, one can easily see that $X^{m}$ has
the following approximate scaling property:
$$
\{m^{-1/\alpha}(X^{1}_{mt}-X_{0}^{1}),t\geq 0\} \text{ has the same law as } \{X_{t}^{m}-X_{0}^{m},t\geq 0\}.
$$
In terms of transition densities, this approximate scaling property can be written as
\beq\label{RSP:scaling1}
p^{m}(t,x,y)=m^{d/\alpha}p^{1}(mt,m^{1/\alpha}x,m^{1/\alpha}y).
\eeq

It is well known that the transition density $p_{D}^{m}(t,x,y)$ of $X^{m, D}$ is continuous
on $(0,\infty)\times D \times D$.
Since both $p^{m}(t,x,y)$ and $p_{D}^{m}(t,x,y)$ are continuous on $(0,\infty)\times D \times D$, $r_{D}^{m}(t,x,y)=p^{m}(t,x,y)-p_{D}^{m}(t,x,y)$
is also continuous there.
$p_{D}^{m}(t,x,y)$ and $r_{D}^{m}(t,x,y)$ also enjoy the following approximate scaling property:
\beq\label{RSP:scaling3}
p^1_{m^{1/\alpha}D}(t,x,y)=m^{-d/\alpha}p_{D}^{m}(t/m,x/m^{1/\alpha},y/m^{1/\alpha}),
\eeq
\beq\label{RSP:scaling5}
r^1_{m^{1/\alpha}D}(t,x,y)=m^{-d/\alpha}r_{D}^{m}(t/m,x/m^{1/\alpha},y/m^{1/\alpha}).
\eeq

The L\'evy measure of the relativistic $\alpha$-stable process $X^{m}$ has a density
$$
J^{m}(x)=j^{m}(|x|):=\frac{\alpha}{2\Gamma(1-\alpha/2)}\int_{0}^{\infty}
(4\pi u)^{-d/2}e^{-|x|^{2}/4u}e^{-m^{2/\alpha}u}u^{-(1+\alpha/2)}du,
$$
which is continuous and radially decreasing on $\R^{d}\setminus \{0\}$ (see \cite[Lemma 2]{Ry}).
Put $J^{m}(x,y):=j^{m}(|x-y|)$.
Let $\mathcal{A}(d,-\alpha):=\alpha2^{\alpha-1}\pi^{-d/2}\Gamma(\frac{d+\alpha}{2})\Gamma(1-\frac{\alpha}{2})^{-1}$. Using change of variables twice, first
with $u=|x|^{2}v$ then with $v=1/s$, we get
\beq\label{Levy:RSP}
J^{m}(x,y)=\mathcal{A}(d,-\alpha)|x-y|^{-d-\alpha}\psi(m^{1/\alpha}|x-y|),
\eeq
where
\beq\label{psi}
\psi(r):=2^{-(d+\alpha)}\Gamma\left(\frac{d+\alpha}{2}\right)^{-1}\int_{0}^{\infty}s^{(d+\alpha)/2-1}e^{-s/4-r^{2}/s}ds,
\eeq
which satisfies $\psi(0)=1$ and
$$
c_{1}^{-1}e^{-r}r^{(d+\alpha-1)/2}\leq\psi(r)\leq c_{1}e^{-r}r^{(d+\alpha-1)/2} \quad \text{on } [1,\infty)
$$
for some $c_{1}>1$ (see \cite[pp. 276-277]{CS3} for details). We denote the L\'evy density of $X$ by
$$
J(x,y):=J^{0}(x,y)=\mathcal{A}(d,-\alpha)|x-y|^{-d-\alpha}.
$$
Note that from \eqref{Levy:RSP} and \eqref{psi} we see that for any $x \in \R^{d}\setminus \{0\}$
$$
j^{m}(|x|)\leq j^{0}(|x|).
$$
It follows from \cite[Theorem 4.1.]{CKS2} that, for any positive constants $M$ and
$T$ there exists a constant $c>1$ such that for all
$m\in(0,M]$, $t\in(0,T]$, and $x,y \in \R^{d}$ we have
\beq\label{RSP:HK}
c^{-1}\left(t^{-d/\alpha}\wedge tJ^{m}(x,y)\right) \leq p^{m}(t,x,y) \leq c\left(t^{-d/\alpha}\wedge tJ^{m}(x,y)\right).
\eeq

We will need a simple lemma from \cite{GR} about the relationship between
$r_{D}^{m}(t,x,y)$ and $r_{D}^{0}(t,x,y)$.
The lemma is true in
much more general situations but we just need it when one of the processes is a symmetric
$\alpha$-stable process and the other is a relativistic $\alpha$-stable process.

\begin{lemma}\label{r:relation}
Suppose that $X$ and $Y$ are two L\'evy processes with L\'evy densities $J^{X}$ and $J^{Y}$, respectively.
Suppose that $\sigma=J^{X}-J^{Y}$ is nonnegative on $\R^d$ with $\int_{\R^d}\sigma(x)dx=\ell<\infty$ and
$D$ is an open set. Then for any $x\in D$ and $t>0$,
$$
p_{D}^{Y}(t,x,\cdot)\leq e^{\ell t}p_{D}^{X}(t,x,\cdot) \quad a.s.
$$
If, in addition, $p^{X}(t,\cdot)$ and $p^{Y}(t,\cdot)$ are continuous, then we have for $x,y \in D$,
$$
r_{D}^{Y}(t,x,y)\leq e^{2\ell t}r_{D}^{X}(t,x,y).
$$
\end{lemma}

The next proposition is the (generalized) Ikeda-Watanabe formula for the relativistic stable process, which describes the joint distribution of
$\tau^{m}_{D}$ and $X^{m}_{\tau_{D}^{m}}$.
\begin{prop}[Proposition 2.7 \cite{KulSi}]\label{IW:RSP}
Assume that $D$ is an open subset of $\R^{d}$ and $A$ is a Borel set
such that $A\subset D^{c}\setminus \partial D$.
If $0\leq t_{1}<t_{2}<\infty$, then
$$
\P_{x}\left(X^{m}_{\tau_{D}^{m}}\in A, \ t_{1}<\tau_{D}^{m}<t_{2}\right)=
\int_{D}\int_{t_{1}}^{t_{2}}p_{D}^{m}(s,x,y)ds\int_{A}J^{m}(y,z)dzdy, \quad x\in D.
$$
\end{prop}

Now we state a simple lemma about the upper bound of $r_{D}^{m}(t,x,y)$,
which is an analogue of \cite[Lemma 2.1]{BK} for stable processes.

\begin{lemma}\label{ub:HK:RSP}
Let $M,T$ be positive constants. Then there exists a constant $c=c(d,\alpha,M,T)$ such that for all $m\in(0,M]$ and $t\in(0,T]$ we have
$$
r_{D}^{m}(t,x,y)\leq c \left(t^{-d/\alpha}\wedge \frac{t\psi(m^{1/\alpha}\delta_{D}(x))}{\delta_{D}(x)^{d+\alpha}}\right).
$$
\end{lemma}
\pf
Since $\psi$ is eventually decreasing and $\psi(0)=1>0$, there exists a constant $c_{1}>0$
such that $\psi(x)\leq c_{1}\psi(y)$ for all $0\leq y \leq x$.
Now from the definition of $r_{D}^{m}(t,x,y)$ and \eqref{RSP:HK} we have
\begin{eqnarray*}
r_{D}^{m}(t,x,y)&=& r_{D}^{m}(t,y,x)\\
&\leq& \E_{y}\left[t>\tau^{m}_{D} ; p^{m}(t-\tau^{m}_{D},X_{\tau^{m}_{D}}^{m},x)\right]\\
&\leq& \E_{y}\left[c \left(t^{-d/\alpha}\wedge \frac{t\psi(m^{1/\alpha}|x-X^{m}_{\tau^{m}_{D}}|)}{|x-X^{m}_{\tau^{m}_{D}}|^{d+\alpha}}\right)\right]\\
&\leq&c c_{1}\left(t^{-d/\alpha}\wedge \frac{t\psi(m^{1/\alpha}\delta_{D}(x))}{\delta_{D}(x)^{d+\alpha}}\right).
\end{eqnarray*}
\qed

We will need two results from \cite{BK}.
The first result is about the difference $p_{F}^{m}(t,x,y)-p_{D}^{m}(t,x,y)$
when $D\subset F$.
The proof in \cite{BK}, given for stable processes,
mainly uses the strong Markov property and
it works for all strong Markov processes with transition densities.

\begin{prop}[Proposition 2.3 \cite{BK}]\label{diff:hks}
 Let $D$ and $F$ be open sets in $\R^{d}$ such that $D\subset F$. Then for any $x,y \in \R^{d}$ we have
$$
p_{F}^{m}(t,x,y)-p_{D}^{m}(t,x,y)=\E_{x}\left[\tau_{D}^{m}<t, X^{m}_{\tau_{D}^{m}}\in F\setminus D :
p_{F}^{m}(t-\tau_{D}^{m},X^{m}_{\tau_{D}^{m}},y)\right].
$$
\end{prop}

Now we introduce some notation.
Recall that if $D$ is a $C^{1,1}$ open set with characteristics $(r_0, \Lambda_0)$,
then for every $x\in D$ with $\delta_{\partial D}(x)<r_{0}$ and $y\in \R^{d}\setminus \bar{D}$
with $\delta_{\partial D}(y)<r_{0}$, there are
$z_{x},z_{y}\in\partial D$ so that $|x-z_{x}|=\delta_{\partial D}(x)$,
$|y-z_{y}|=\delta_{\partial D}(y)$ and that $B(x_{0},r_{0})\subset D$
and $B(y_{0},r_{0})\subset \R^{d}\setminus\bar{D}$, where $x_{0}=z_{x}+r_{0}(x-z_{x})/|x-z_{x}|$
and $y_{0}=z_{y}+r_{0}(y-z_{y})/|y-z_{y}|$.
Let $H(x)$ be the half-space containing $B(x_{0},r_{0})$ such that $\partial H(x)$ contains
$z_{x}$ and is perpendicular to the segment
$\overline{z_{x}z_{y}}$.
The next proposition says that, in case of the symmetric $\alpha$-stable process,
for small $t$, the quantity $r_{D}^{0}(t,x,x)$ can be replaced by
$r_{H(x)}^{0}(t,x,x)$, which was a very crucial step in proving the main result in \cite{BK}.

\begin{prop}[Proposition 3.1 of \cite{BK}]\label{rDrH:SSP}
Let $D\subset \R^{d}$ be a $C^{1,1}$ open set with characteristics $(r_0, \Lambda_0)$.
Then, for any $x$ with $\delta_{\partial D}(x)<r_{0}/2$ and $t>0$ with
$t^{1/\alpha}\leq r_{0}/2$, we have
$$
\left|r_{D}^{0}(t,x,x)-r_{H(x)}^{0}(t,x,x)\right|\leq \frac{ct^{1/\alpha}}{r_{0}t^{d/\alpha}}
\left(\left(\frac{t^{1/\alpha}}{\delta_{\partial D}(x)}\right)^{d+\frac{\alpha}{2}-1} \wedge 1\right).
$$
\end{prop}

We will need some facts about the ``stability'' of the surface area of the boundary of $C^{1,1}$ open sets.
The following lemma is \cite[Lemma 5]{vB}.

\begin{lemma}
 Let $D$ be a bounded $C^{1,1}$ open set in $\R^{d}$ with characteristic $(r_0, \Lambda_0)$ and define for $0\leq q <r_{0}$,
$$
D_{q}=\{x\in D : \delta_{D}(x)>q\}.
$$
Then
$$
\left(\frac{r_{0}-q}{r_{0}}\right)^{d-1}|\partial D|\leq |\partial D_{q}| \leq \left(\frac{r_{0}}{r_{0}-q}\right)^{d-1}|\partial D|,
\quad 0\leq q < r_{0}.
$$
\end{lemma}

The following result is \cite[Corollary 2.14]{BK}.

\begin{lemma}\label{pertub:C11}
Let $D$ be a bounded $C^{1,1}$ open set in $\R^{d}$ with characteristic $(r_0, \Lambda_0)$.
For any $0<q\leq r_{0}/2$, we have
\begin{enumerate}
 \item $2^{-d+1}|\partial D|\leq |\partial D_{q}|\leq 2^{d-1}|\partial D|$,
\item $|\partial D|\leq \frac{2^{d}|D|}{r_{0}}$,
\item $\left| |\partial D_{q}|-|\partial D|\right|\leq \frac{2^{d}dq|\partial D|}{r_{0}}\leq \frac{2^{2d}dq|D|}{r_{0}^{2}}$.
\end{enumerate}
\end{lemma}

\section{Proof of Theorem \ref{RSP:main}}\label{Proof of the main result}

We first prove that $\lim_{t\rightarrow 0}t^{\frac{d}{\alpha}}Z_{D}^{m}(t)$
exists and identify the limit.

\begin{lemma}\label{RSP,c1 2}
The limit $\lim_{t\rightarrow 0}t^{\frac{d}{\alpha}}Z_{D}^{m}(t)$ exists
and is equal to $C_1|D|$, where $C_1$ is the constant in Theorem \ref{RSP:main}.
\end{lemma}
\pf By definition,
\begin{eqnarray}\label{RSP,c1 3}
&&t^{d/\alpha}Z_{D}^{m}(t)=t^{d/\alpha}\int_{D}p_{D}^{m}(t,x,x)dx\nonumber\\
&=&t^{d/\alpha}\left(\int_{D}p^{m}(t,x,x)dx-\int_{D}r^{m}_{D}(t,x,x)dx\right).
\end{eqnarray}
For the first integral on the right hand side of \eqref{RSP,c1 3}, note that,
by the approximate scaling property
\eqref{RSP:scaling3} and the dominated convergence theorem,
we have, as $t\to 0$,
\begin{eqnarray*}
&&t^{d/\alpha}\left(\int_{D}p^{m}(t,x,x)dx\right)=\int_{D}p^{tm}(1,x,x)dx=|D|p^{tm}(1,0)\\
&\rightarrow&|D|\cdot p^{0}(1,0)=|D|\cdot \frac{\Gamma(d/\alpha) \omega_{d}}{(2\pi)^{d}\alpha}.
\end{eqnarray*}

It remains to show that $\lim_{t\to 0}t^{d/\alpha}\int_{D}r_{D}^{m}(t,x,x)dx=0$.
By Lemma \ref{ub:HK:RSP} we have that
$$
t^{d/\alpha}r_{D}^{m}(t,x,y)\leq c, \qquad (t, x, y)\in(0, 1]\times D\times D,
$$
for some $c>0$.
Hence we have by the monotone convergence theorem,
$$
t^{d/\alpha}\int_{D\setminus D_{t^{1/2\alpha}}}r_{D}^{m}(t,x,x)\rightarrow 0
\quad \text{ as } t \rightarrow 0.
$$
For $x\in D_{t^{1/2\alpha}}$ we have by Lemma \ref{ub:HK:RSP} again for $t\in (0, 1]$,
$$
r_{D}^{m}(t,x,x)\leq c\,t^{\frac12+\frac{d}{2\alpha}}, \qquad x\in D_{t^{1/2\alpha}}.
$$
Hence $\lim_{t\to 0}t^{d/\alpha}\int_{D_{t^{1/2\alpha}}}r_{D}^{m}(t,x,x)dx = 0$.
\qed

It follows from Lemma \ref{RSP,c1 2} that if $N^{m}(\lambda)$ denotes the number of eigenvalues $\lambda^{(m)}_j$ such that $\lambda^{m}_j\le \lambda$,
then it follows from the
classical Karamata Tauberian theorem (see for example \cite{F}) that
$$
N^{m}(\lambda) \sim \frac{C_{1}|D|}{\Gamma(d/\alpha +1)}\lambda^{d/\alpha}, \quad \text{as } \lambda \rightarrow \infty.
$$
This is the analogue for killed relativistic stable processes of the celebrated Weyl's asymptotic formula for
the eigenvalues of the Dirichlet Laplacian and it is already proved in \cite{BMN} (see \cite[(1.10)]{BMN}).
This result has been known at least since 2009, see \cite[Remark 1.2]{BMN}.

Now we focus on identifying the next terms in $Z_{D}^{m}(t)$.
For this, we need to find the order of $t$ in
$Z_{D}^{m}(t)-C_{1}t^{-\frac{d}{\alpha}}$.
Note that by Lemma \ref{RSP,c1 2},
\begin{eqnarray*}
&&Z_{D}^{m}(t)-C_{1}t^{-d/\alpha}=\int_{D}p_{D}^{m}(t,x,x)-p^{0}(t,x,x)dx\\
&=&\int_{D}\left(p^{m}(t,x,x)-p^{0}(t,x,x)\right)dx-\int_{D}r_{D}^{m}(t,x,x)dx.
\end{eqnarray*}

The next lemma gives the orders of $t$ in
$p^{m}(t,x,x)-p^{0}(t,x,x)$ up to $t^{\frac{2}{\alpha}}t^{-\frac{d}{\alpha}}$.

\begin{lemma}\label{RSP:power terms}
Let $k$ be the largest integer such that $k<\frac{2}{\alpha}$. Then we have
$$
p^{m}(t,x,x)-p^{0}(t,x,x)=t^{-d/\alpha}\frac{\omega_{d}\Gamma(d/\alpha)}{(2\pi)^{d}\alpha}\sum_{n=1}^{k}\frac{m^{n}}{n!}t^{n}+O(t^{2/\alpha}t^{-d/\alpha}).
$$
\end{lemma}
\pf
By the scaling property \eqref{RSP:scaling1} we have
\begin{eqnarray*}
&&p^{m}(t,x,x)-p^{0}(t,x,x)=p^{m}(t,0)-p^{0}(t,0)\\
&=&t^{-d/\alpha}\left(p^{tm}(1,0)-p^{0}(1,0)\right)\\
&=&t^{-d/\alpha}(2\pi)^{-d}\int_{\R^{d}}e^{-(|\xi|^{2}+(mt)^{2/\alpha})^{\alpha/2}+mt}-e^{-|\xi|^{\alpha}}d\xi.
\end{eqnarray*}
Note that for any $x\geq 0$ we have $(1+x)^{\alpha/2}\leq 1 +\frac{\alpha}2 x$. Thus
\begin{eqnarray*}
&&\left(|\xi|^{2}+(mt)^{2/\alpha}\right)^{\alpha/2}=|\xi|^{\alpha}
\left(1+\frac{(mt)^{2/\alpha}}{|\xi|^{2}}\right)^{\alpha/2}
\leq |\xi|^{\alpha}\left(1+\frac{\alpha}2\frac{(mt)^{2/\alpha}}{|\xi|^{2}}\right).\nonumber\\
\end{eqnarray*}
Consequently
\begin{eqnarray*}
&&0\leq e^{-|\xi|^{\alpha}}-e^{-\left(|\xi|^{2}+(mt)^{2/\alpha}\right)^{\alpha/2}}\nonumber\\
&\leq& e^{-|\xi|^{\alpha}}-e^{-|\xi|^{\alpha}\left(1+\frac{\alpha}2\frac{(mt)^{2/\alpha}}{|\xi|^{2}}\right)}=
e^{-|\xi|^{\alpha}}\left(1-e^{-\frac{\alpha}2\frac{(mt)^{2/\alpha}}{|\xi|^{2-\alpha}}}\right)\nonumber\\
&\leq&e^{-|\xi|^{\alpha}}\left(\frac{\alpha}2\frac{(mt)^{2/\alpha}}{|\xi|^{2-\alpha}}\right),
\end{eqnarray*}
where we used $1-e^{-x}\leq x$ for all $x\geq 0$ in the last inequality above.
Therefore
\begin{eqnarray*}
&&0\leq \int_{\R^{d}}e^{-\left(|\xi|^{2}+(mt)^{2/\alpha}\right)^{\alpha/2}+mt}-e^{-|\xi|^{\alpha}}d\xi\\
&\leq&\int_{\R^{d}}\left|e^{-\left(|\xi|^{2}+(mt)^{2/\alpha}\right)^{\alpha/2}+mt}-e^{-|\xi|^{\alpha}}e^{mt}
+e^{-|\xi|^{\alpha}}e^{mt}-e^{-|\xi|^{\alpha}}\right|d\xi\\
&\leq&\int_{\R^{d}}\left|e^{-\left(|\xi|^{2}+(mt)^{2/\alpha}\right)^{\alpha/2}+mt}-e^{-|\xi|^{\alpha}}e^{mt}\right|d\xi
+\int_{\R^{d}}\left|e^{-|\xi|^{\alpha}}e^{mt}-e^{-|\xi|^{\alpha}}\right|d\xi\\
&\leq&\int_{\R^{d}}e^{mt}e^{-|\xi|^{\alpha}}\left(\frac{\alpha}2\frac{(mt)^{2/\alpha}}{|\xi|^{2-\alpha}}\right)d\xi
+\int_{\R^{d}}e^{-|\xi|^{\alpha}}\left(e^{mt}-1\right)d\xi\\
&=&e^{mt}\frac{\alpha}2(mt)^{2/\alpha}\int_{\R^{d}}\frac{e^{-|\xi|^{\alpha}}}{|\xi|^{2-\alpha}}d\xi
+\sum_{n=1}^{\infty}\frac{(mt)^{n}}{n!}\int_{\R^{d}}e^{-|\xi|^{\alpha}}d\xi.
\end{eqnarray*}
Since $k+j \geq 2/\alpha$ for any $j\geq 1$, we have $\displaystyle\sum_{n=k+1}^{\infty}\frac{(mt)^{n}}{n!}=O(t^{2/\alpha})$.
Therefore
$$
\int_{\R^{d}}\left(e^{-\left(|\xi|^{2}+(mt)^{2/\alpha}\right)^{\alpha/2}+mt}-e^{-|\xi|^{\alpha}}\right)d\xi
=O(t^{2/\alpha})+\frac{\omega_{d}\Gamma(d/\alpha)}{\alpha}\sum_{n=1}^{k}\frac{(mt)^{n}}{n!}
$$
and
$$
p^{m}(t,x,x)-p^{0}(t,x,x)=t^{-d/\alpha}\frac{\omega_{d}\Gamma(d/\alpha)}{(2\pi)^{d}\alpha}\sum_{n=1}^{k}\frac{m^{n}}{n!}t^{n}+O(t^{2/\alpha}t^{-d/\alpha}).
$$
\qed

Now we try to find the orders of $t$ in the expansion of $\int_{D}r_{D}^{m}(t,x,x)dx$ up
to the order of $t^{\frac{2}{\alpha}}t^{-\frac{d}{\alpha}}$.
For this,
we need to assume some regularity condition on the boundary of $D$.
Hence in the remainder of this section we assume that $D$ is a bounded $C^{1,1}$ open set with characteristic $(r_0, \Lambda_0)$.
We also assume that $t^{1/\alpha}\leq\frac{r_{0}}{2}$.

We first deal with the contribution in $D_{r_{0}/2}$.
\begin{lemma}\label{RSP:2nd term 1}
There exists $c=c(d, \alpha)>0$ such that
$$
\dint_{D_{r_{0}/2}} r_{D}^{m}(t,x,x)dx\leq c e^{2mt}\frac{|D|t^{2/\alpha}}{r_{0}^{2}t^{d/\alpha}}.
$$
\end{lemma}
\pf
It follows from Lemma \ref{r:relation} that $r_{D}^{m}(t,x,y)\leq e^{2mt}r_{D}^{0}(t,x,y)$.
By \cite[(3.2)]{BK}
we know that
$$
\int_{D_{r_{0}/2}}r_{D}^{0}(t,x,y)dx\leq \frac{c|D|t^{2/\alpha}}{r_{0}^{2}t^{d/\alpha}}.
$$
The desired assertion follows immediately.
\qed
\begin{lemma}\label{RSP:2nd term 2}
There exists $c=c(d, \alpha)>0$ such that
$$
r_{D}^{m}(t,x,x)-r_{H(x)}^{m}(t,x,x)\leq c e^{2mt}\frac{t^{1/\alpha}}{t^{d/\alpha}}\left(\left(\frac{t^{1/\alpha}}{\delta_{D}(x)}
 \right)^{d+\frac{\alpha}{2}-1}\wedge 1\right)
$$
and
$$\dint_{D\setminus D_{r_{0}/2}}\left(r_{D}^{m}(t,x,x)-r_{H(x)}^{m}(t,x,x)\right)dx \leq ce^{2mt}\frac{t^{2/\alpha}}{t^{d/\alpha}}.
$$
\end{lemma}
\pf
If the first assertion of the lemma is right, then it is easy to see that
$$
\int_{D\setminus D_{r_{0}/2}}\left(\left(\frac{t^{1/\alpha}}{\delta_{D}(x)}\right)^{d+\frac{\alpha}{2}-1}\wedge 1\right)dx\leq ct^{1/\alpha}.
$$
Hence we focus on proving the first assertion. By \cite[(3.4)]{BK}, we know that
$$
r_{D}^{0}(t,x,x)-r_{H(x)}^{0}(t,x,x)\leq c\frac{t^{1/\alpha}}{t^{d/\alpha}}\left(\left(\frac{t^{1/\alpha}}{\delta_{D}(x)}\right)^{d+\frac{\alpha}{2}-1}\wedge 1\right).
$$
Recall that $J^{m}(x)\leq J^{0}(x)$ for any $x\in \R^{d}\setminus \{0\}$.
Now it follows from the generalized Ikeda-Watanabe formula and Lemma \ref{r:relation} that
\begin{eqnarray*}
&&r_{D}^{m}(t,x,x)-r_{H(x)}^{m}(t,x,x)\\
&=&\E_{x}\left[t>\tau_{D}^{m}, X^{m}_{\tau^{m}_{D}}\in H(x)\setminus D; p_{H(x)}^{m}(t-\tau^{m}_{D}, X^{m}_{\tau^{m}_{D}},x)\right]\\
&=&\int_{D}\int_{0}^{t}p_{D}^{m}(s,x,y)ds\int_{H(x)\setminus D}J^{m}(y,z)p_{H(x)}^{m}(t-s,z,x)dzdy\\
&\leq& e^{2mt}\int_{D}\int_{0}^{t}p_{D}^{0}(s,x,y)ds\int_{H(x)\setminus D}J^0(y,z)p_{H(x)}^{0}(t-s,z,x)dzdy\\
&=&e^{2mt}\E_{x}\left[t>\tau_{D}^{0}, X_{\tau^{0}_{D}}\in H(x)\setminus D; p_{H(x)}^{0}(t-\tau^{0}_{D}, X_{\tau^{0}_{D}},x)\right]\\
&=&e^{2mt}\left(r_{D}^{0}(t,x,x)-r_{H(x)}^{0}(t,x,x)\right)\\
&\leq&ce^{2mt}\frac{t^{1/\alpha}}{t^{d/\alpha}}\left((\frac{t^{1/\alpha}}{\delta_{D}(x)})^{d+\frac{\alpha}{2}-1}\wedge 1\right).
\end{eqnarray*}
\qed
\begin{lemma}\label{RSP:2nd term 3}
There exists $c=c(d, \alpha)>0$ such that
$$
\dint_{D\setminus D_{r_{0}/2}}r_{H(x)}^{m}(t,x,x)dx-t^{1/\alpha}t^{-d/\alpha}\dint_{0}^{\frac{r_{0}}{2t^{1/\alpha}}}|\partial D|f_{H}^{tm}(1,u)du
\leq c t^{2/\alpha}t^{-d/\alpha}.
$$
\end{lemma}
\pf
Using the scaling relation \eqref{RSP:scaling5} we get
\begin{eqnarray*}
&&\int_{D\setminus D_{r_{0}/2}}r_{H(x)}^{m}(t,x,x)dx\\
&=&\int_{0}^{r_{0}/2}|\partial D_{u}|f_{H}^{m}(t,u)du\\
&=& \int_{0}^{r_{0}/2}|\partial D_{u}|t^{-d/\alpha}f_{H}^{tm}(1,u/t^{1/\alpha})du\\
&=&t^{1/\alpha}t^{-d/\alpha}\int_{0}^{r_{0}/2t^{1/\alpha}}|\partial D_{ut^{1/\alpha}}|f_{H}^{tm}(1,u)du.
\end{eqnarray*}
It follows from Corollary \ref{pertub:C11} that $\left||\partial D_{q}|-|\partial D|\right|\leq \frac{2^{d}dq|\partial D|}{r_{0}}\leq \frac{2^{2d}dq|D|}{r_{0}^{2}}$ for
any $q\leq r_{0}/2$. Hence
\begin{eqnarray*}
&&\left|\int_{D\setminus D_{r_{0}/2}}r_{H(x)}^{m}(t,x,x)-t^{1/\alpha}t^{-d/\alpha}\int_{0}^{\frac{r_{0}}{2t^{1/\alpha}}}|\partial D|
f_{H}^{tm}(1,u)du\right|\\
&\leq&t^{1/\alpha}t^{-d/\alpha}\int_{0}^{\frac{r_{0}}{2t^{1/\alpha}}}\left||\partial D_{ut^{1/\alpha}}|-|\partial D|\right|f_{H}^{tm}(1,u)du\\
&\leq&c_{1} t^{2/\alpha}t^{-d/\alpha}\int_{0}^{\infty}uf_{H}^{tm}(1,u)du\\
&\leq&c_{2} t^{2/\alpha}t^{-d/\alpha}.
\end{eqnarray*}
\qed

\begin{lemma}\label{RSP:2nd term 4}
There exists $c=c(d, \alpha)>0$ such that
$$
t^{1/\alpha}t^{-d/\alpha}\dint_{0}^{\infty}|\partial D|f_{H}^{tm}(1,u)du -
 t^{1/\alpha}t^{-d/\alpha}\dint_{0}^{\frac{r_{0}}{2t^{1/\alpha}}}|\partial D|f_{H}^{tm}(1,u)du \leq ct^{2/\alpha}t^{-d/\alpha}.
$$
\end{lemma}
\pf
It follows from Lemma \ref{r:relation} that
\begin{eqnarray*}
&&t^{1/\alpha}t^{-d/\alpha}\dint_{0}^{\infty}|\partial D|f_{H}^{tm}(1,u)du -
t^{1/\alpha}t^{-d/\alpha}\dint_{0}^{\frac{r_{0}}{2t^{1/\alpha}}}|\partial D|f_{H}^{tm}(1,u)du \\
&=&t^{1/\alpha}t^{-d/\alpha}\dint_{\frac{r_{0}}{2t^{1/\alpha}}}^{\infty}|\partial D|f_{H}^{tm}(1,u)du\\
&=&t^{1/\alpha}t^{-d/\alpha}|\partial D|\dint_{\frac{r_{0}}{2t^{1/\alpha}}}^{\infty}f_{H}^{tm}(1,u)du\\
&\leq&e^{2mt}t^{1/\alpha}t^{-d/\alpha}|\partial D|\dint_{\frac{r_{0}}{2t^{1/\alpha}}}^{\infty}f_{H}^{0}(1,u)du.
\end{eqnarray*}
For $q\geq r_{0}/(2t^{1/\alpha})$ we have $f_{H}^{0}(1,q)\leq c q^{-d-\alpha}\leq cq^{-2}$. Hence
$$
\dint_{\frac{r_{0}}{2t^{1/\alpha}}}^{\infty}f_{H}^{0}(1,u)du\leq c\dint_{\frac{r_{0}}{2t^{1/\alpha}}}^{\infty}\frac{dq}{q^{2}}\leq c\frac{t^{1/\alpha}}{r_{0}}
$$
and the result now follows.
\qed

\begin{lemma}\label{RSP:2nd term 5}
 $\displaystyle\lim_{t \downarrow 0}\int_{0}^{\infty}f_{H}^{tm}(1,u)du=\int_{0}^{\infty}f_{H}^{0}(1,u)du$.
\end{lemma}
\pf
This follows immediately from the continuity of $m \mapsto r_{D}^{m}(t,x,y)$ and the dominated convergence theorem.
\qed

{\bf Proof of Theorem \ref{RSP:main}}
Combining Lemmas \ref{RSP,c1 2}, \ref{RSP:power terms}, \ref{RSP:2nd term 1}, \ref{RSP:2nd term 2},
\ref{RSP:2nd term 3}, \ref{RSP:2nd term 4}, and \ref{RSP:2nd term 5}, we immediately arrive at
Theorem \ref{RSP:main}.
\qed

\section{Proof of Theorem \ref{RSP:main2}}\label{main result 2:Lipschitz}

In this section we always assume that $D$ is a bounded Lipschitz open set in $\R^d$.
The argument of this section is similar to previous section and \cite{BKS}.
We will follow the argument in \cite{BKS} closely,
making necessary modifications for relativistic stable processes.
Note that even though the main theorem in \cite{BKS} is stated for a Lipschitz domain,
it remains true for a bounded Lipschitz open set.

First we need two technical facts which play crucial roles later.
The first proposition is \cite[Proposition 2.9]{BKS} and
we will state it here for reader's convenience.

\begin{prop}[Proposition 2.9. \cite{BKS}]\label{conv1}
Suppose that $f:(0,\infty)\rightarrow [0, \infty)$ is continuous and satisfies
$f(r)\leq c(1\wedge r^{-\beta})$ for some $\beta>1$.
Furthermore, suppose that for any $0<R_{1}<R_{2}<\infty$, $f$ is
Lipschitz on $[R_{1},R_{2}]$. Then we have
$$
\lim_{\eta\rightarrow 0^{+}}\frac{1}{\eta}\int_{D}f\left(\frac{\delta_{D}(x)}{\eta}\right)dx=\mathcal{H}^{d-1}(\partial D)\int_{0}^{\infty}f(r)dr.
$$
\end{prop}

\begin{lemma}\label{conv2}
Suppose that $f:(0,\infty)\rightarrow [0, \infty)$ is continuous and satisfies
$f(r)\leq c_1(1\wedge r^{-\beta})$ for some $\beta>1$.
Furthermore, suppose that for any $0<R_{1}<R_{2}<\infty$, $f$ is
Lipschitz on $[R_{1},R_{2}]$.
Let $\{f^{\eta}:\eta>0\}$ be continuous functions from $(0, \infty)$ to $[0, \infty)$
such that, for any $0<L<M<\infty$,
$\displaystyle\lim_{\eta\rightarrow 0}f^{\eta}(r)=f(r)$ uniformly for $r\in [L,M]$.
Suppose that there exists $c_2>0$ such that
$f^{\eta}(r)\leq c_2f(r)$ for all $\eta\leq 1$. Then we have
$$
\lim_{\eta\rightarrow 0^{+}}\frac{1}{\eta}\int_{D}f^{\eta}\left(\frac{\delta_{D}(x)}{\eta}\right)dx=\mathcal{H}^{d-1}(\partial D)\int_{0}^{\infty}f(r)dr.
$$
\end{lemma}
\pf
Let $\psi_{\eta}(r)=\eta^{-1}\left|\{x\in D : \delta_{D}(x)<\eta r\}\right|$.
Note (cf. proof of \cite[Proposition 1.1]{B}) that $\psi_\eta(r)\le c$ for all $\eta, r>0$
and that
$$
\eta^{-1}\int_{D}f\left(\frac{\delta_{D}(x)}{\eta}\right)dx=\int_{0}^{\infty}f(r)d\psi_{\eta}(r),\nonumber
$$
and
$$
\eta^{-1}\int_{D}f^{\eta}\left(\frac{\delta_{D}(x)}{\eta}\right)dx=\int_{0}^{\infty}f^{\eta}(r)d\psi_{\eta}(r).\nonumber
$$
It was shown in \cite[Proposition 2.9.]{BKS} that, for any $0<R_{1}<R_{2}<\infty$ and $\eta>0$, $f$ satisfies
\beq\label{f:2}
\int_{0}^{R_{1}}f(r)d\psi_{\eta}(r)\leq cR_{1},
\eeq
\beq\label{f:3}
\int_{R_{2}}^{\infty}f(r)d\psi_{\eta}(r)\leq c\eta^{\beta-1}+cR_{2}^{1-\beta},
\eeq
$$
\lim_{\eta\rightarrow 0^{+}}\int_{R_{1}}^{R_{2}}f(r)d\psi_{\eta}(r)=\mathcal{H}^{d-1}(\partial D)\int_{R_{1}}^{R_{2}}f(r)dr.
$$
Since $f^{\eta}\leq c_2f$ for $\eta \leq 1$ we have the same inequalities as \eqref{f:2} and \eqref{f:3} for $f^{\eta}$, $\eta\leq 1$.
Hence it is enough to show that
$$
\lim_{\eta\rightarrow 0^{+}}\int_{R_{1}}^{R_{2}}f^{\eta}(r)d\psi_{\eta}(r)=\mathcal{H}^{d-1}(\partial D)\int_{R_{1}}^{R_{2}}f(r)dr.
$$
For any partition $R_{1}=x_{0}<x_{1}<\cdots<x_{n}=R_{2}$ of $[R_{1},R_{2}]$, we have
\begin{eqnarray*}
 &&\left|\sum_{i=1}^{n}f^{\eta}(x_{i})\left(\psi_{n}(x_{i})-\psi_{n}(x_{i-1})\right)-
\sum_{i=1}^{n}f(x_{i})\left(\psi_{n}(x_{i})-\psi_{n}(x_{i-1})\right)\right|\nonumber\\
&=&\sum_{i=1}^{n}|f^{\eta}(x_{i})-f(x_{i})|\left(\psi_{n}(x_{i})-\psi_{n}(x_{i-1})\right)\nonumber\\
&\leq&\parallel f^{\eta}-f\parallel_{L^{\infty}[R_{1},R_{2}]}\psi_{\eta}(R_{2}).
\end{eqnarray*}
Note that for any $\eta>0$ the function $r\rightarrow \psi_{\eta}(r)$ is nondecreasing and for any $\eta>0$, $r>0$ we have
$\psi_{\eta}(r)\leq cr$ for some constant $c$.
Since $f^{\eta}\rightarrow f$ uniformly on $r\in [R_{1},R_{2}]$, taking supremum for all possible partitions gives
$$
\lim_{\eta\rightarrow 0^{+}}\int_{R_{1}}^{R_{2}}f^{\eta}(r)d\psi_{\eta}(r)=\lim_{\eta\rightarrow 0^{+}}\int_{R_{1}}^{R_{2}}f(r)d\psi_{\eta}(r)
=\mathcal{H}^{d-1}(\partial D)\int_{R_{1}}^{R_{2}}f(r)dr.
$$
\qed

\begin{lemma}\label{conv:unif2}
For any $0<L<M<\infty$,
$p^{m}(t,x,y)$ converges uniformly to $p^{0}(t,x,y)$ as $m\rightarrow 0$ for $(t,x,y)\in [L,M]\times \R^{d}\times \R^{d}$.
\end{lemma}
\pf
Note that
\begin{eqnarray*}
&&\left|p^{m}(t,x,y)-p^{0}(t,x,y)\right|=\left|
(2\pi)^{-d}\int_{\R^{d}}e^{-i\xi(y-x)}\left(e^{-t((|\xi|^{2}+m^{2/\alpha})^{\alpha/2}-m)}-e^{-t|\xi|^{\alpha}}\right)d\xi\right|\\
&\leq&(2\pi)^{-d}\int_{\R^{d}}\left|e^{-i\xi(y-x)}\left(e^{-t((|\xi|^{2}+m^{2/\alpha})^{\alpha/2}-m)}-e^{-t|\xi|^{\alpha}}\right)\right|d\xi\\
&\leq&(2\pi)^{-d}\int_{\R^{d}}e^{-t((|\xi|^{2}+m^{2/\alpha})^{\alpha/2}-m)}-e^{-t|\xi|^{\alpha}}d\xi\\
&=&(2\pi)^{-d}(p^{m}(t,0)-p^0(t,0)).
\end{eqnarray*}
Now it follows from Lemma \ref{RSP:power terms} that for $t\in[L,M]$ and $x,y\in \R^{d}$,
\begin{eqnarray*}
 &&\left|p^0(t,x,y)-p^{m}(t,x,y)\right|\\
&\leq&t^{-d/\alpha}(2\pi)^{-d}e^{mt}\frac{\alpha}2(mt)^{2/\alpha}\int_{\R^{d}}\frac{e^{-|\xi|^{\alpha}}}{|\xi|^{2-\alpha}}d\xi
+\sum_{n=1}^{\infty}\frac{(mt)^{n}}{n!}\int_{\R^{d}}e^{-|\xi|^{\alpha}}d\xi\\
&\leq&L^{-d/\alpha}(2\pi)^{-d}\frac{\alpha}2(mM)^{2/\alpha}\int_{\R^{d}}\frac{e^{-|\xi|^{\alpha}}}{|\xi|^{2-\alpha}}d\xi
+\sum_{n=1}^{\infty}\frac{(mM)^{n}}{n!}\int_{\R^{d}}e^{-|\xi|^{\alpha}}d\xi.
\end{eqnarray*}
The last quantity above converges to $0$ as $m\rightarrow 0$.
\qed

For convenience, we define the following notation.
$$
f_{H}^{m}(t,r):=r_{H}^{m}(t,(r, \tilde{0}), (r, \tilde{0})), \qquad r>0.
$$

\begin{lemma}\label{conv:unif}
For any $0<L<M<\infty$ and $m>0$,
$$
\lim_{t\rightarrow 0}f_{H}^{tm}(1,r)=f^0_{H}(1,r), \quad \text{uniformly in } r\in[L,M],
$$
that is, given $\eps>0$ there exists $t_{0}>0$ such that for $0\leq t\leq t_{0}$ we have
$$
\sup_{r\in[L,M]}\left|r_{H}^{tm}(1,(r,\tilde{0}),(r,\tilde{0}))-r^0_{H}(1,(r,\tilde{0}),(r,\tilde{0}))\right|<\eps.
$$
\end{lemma}
\pf
Recall that $r^0_{H}(t,x,y)=\E_{x}[\tau^0_{H}<t,p^0(t-\tau^0_{H},X_{\tau^0_{H}},y)]$ and
$r_{H}^{m}(t,x,y)=\E_{x}[\tau^{m}_{H}<t,p^{m}(t-\tau^{m}_{H},X^{m}_{\tau^{m}_{H}},y)]$.
It is well known that
$$
p^0(t,x,y)\asymp t^{-d/\alpha}\wedge \frac{t}{|x-y|^{d+\alpha}}.
$$
Since $|X^0_{\tau_{H}}-(r,\tilde{0})|>L$, we have, together with Lemma \ref{r:relation},
$$
p^{0}(1-\tau^0_{H},X^0_{\tau_{H}^{0}},(r,\tilde{0}))\leq c \frac{1-\tau^0_{H}}{L^{d+\alpha}},
$$
$$
p^{tm}(1-\tau^{tm}_{H},X^{tm}_{\tau_{H}^{tm}},(r,\tilde{0}))\leq c e^{tm}\frac{1-\tau_{H}^{tm}}{L^{d+\alpha}}.
$$
Now take $\delta_{1}$ small so that
\beq\label{unif:r1}
\E_{(r,\tilde{0})}\left[1-\delta_{1}\leq \tau^0_{H}<1, p^0(1-\tau^0_{H},X^0_{\tau^0_{H}},(r,\tilde{0}))\right]< \eps,
\eeq
\beq \label{unif:r2}
\E_{(r,\tilde{0})}\left[1-\delta_{1}\leq \tau_{H}^{tm}<1,
p^{tm}(1-\tau^{tm}_{H},X^{tm}_{\tau_{H}^{tm}},(r,\tilde{0}))\right]< \eps.
\eeq

Now let $V^{m}$ be a L\'evy process with L\'evy density $\sigma=J-J^{m}$ and define
$T^{m}:=\inf\{t>0 : V^{m}_{t}\neq 0\}$.
Then $V^{m}$ is a compound Poisson process and $T^{m}$ is an exponential random variable with parameter $m$
and independent of $X$ (See \cite{Ry}).
Then we have
\begin{eqnarray*}
 &&\E_{(r,\tilde{0})}\left[\tau_{H}^{tm}<1-\delta_{1}, p^{tm}(1-\tau_{H}^{tm},X^{tm}_{\tau_{H}^{tm}},(r,\tilde{0}))\right]\\
&=&\E_{(r,\tilde{0})}\left[T^{tm}>1, \tau_{H}^{tm}<1-\delta_{1}, p^{tm}(1-\tau_{H}^{tm},X^{tm}_{\tau_{H}^{tm}},(r,\tilde{0}))\right] \\
&+& \E_{(r,\tilde{0})}\left[T^{tm}\leq 1,\tau_{H}^{tm}<1-\delta_{1}, p^{tm}(1-\tau_{H}^{tm},X^{tm}_{\tau_{H}^{tm}},(r,\tilde{0}))\right].
\end{eqnarray*}
Since $p^{tm}(1-\tau_{H}^{tm},X^{tm}_{\tau_{H}^{tm}},(r,\tilde{0})) \leq c \frac{e^{mt}}{L^{d+\alpha}}$, we have
\beq\label{unif:r3}
\E_{(r,\tilde{0})}\left[T^{tm}\leq 1,\tau_{H}^{tm}<1-\delta_{1}, p^{tm}(1-\tau_{H}^{tm},X^{tm}_{\tau_{H}^{tm}},(r,\tilde{0}))\right]\leq
c\frac{e^{mt}}{L^{d+\alpha}}(1-e^{-mt}).
\eeq
Similarly we also have
\beq\label{unif:r4}
\E_{(r,\tilde{0})}\left[T^{tm}\leq 1,\tau^0_{H}<1-\delta_{1}, p^0(1-\tau^0_{H},X^0_{\tau_{H}},(r,\tilde{0}))\right]\leq
c \frac{1}{L^{d+\alpha}}(1-e^{-mt}).
\eeq

Take $t_{1}>0$ such that \eqref{unif:r3} and \eqref{unif:r4} is less than $\eps$ for all $t\leq t_{1}$.
Next note that for $T^{tm}>1$ and $\tau_{H}^{tm}<1$,
we have $\tau_{H}^{tm}=\tau^0_{H}$ and $X^{tm}_{\tau_{H}^{tm}}=X^0_{\tau^0_{H}}$.
Hence it follows that
\begin{eqnarray}\label{unif:r5}
&& | \E_{(r,\tilde{0})}\left[T^{tm}>1, \tau^{tm}_{H}<1-\delta_{1}, p^{tm}(1-\tau^{tm}_{H},X^{tm}_{\tau_{H}^{tm}},(r,\tilde{0}))\right] \nonumber\\
&&- \E_{(r,\tilde{0})}\left[T^{tm}>1,\tau^0_{H}<1-\delta_{1}, p^0(1-\tau^0_{H},X^0_{\tau^0_{H}},(r,\tilde{0}))\right]| \nonumber\\
&\leq&\E_{(r,\tilde{0})}\left[ T^{tm}>1, \tau^0_{H}<1-\delta_{1}, |p^{tm}(1-\tau^0_{H},X^0_{\tau^0_{H}},(r,\tilde{0}))
-p^0(1-\tau^0_{H},X^0_{\tau^0_{H}},(r,\tilde{0}))| \right] \nonumber\\
&\leq& \sup_{s\in[\delta_{1},1], x,y \in \R^{d}}|p^{tm}(s,x,y)-p^0(s,x,y)|.
\end{eqnarray}
It follows from Lemma \ref{conv:unif2} that there exists $t_{2}>0$ such that $\sup_{s\in[\delta_{1},1],
x,y \in \R^{d}}|p^{tm}(s,x,y)-p^0(s,x,y)|<\eps$ for $0\leq t\leq t_{2}$.
Now let $t_{0}=t_{1}\wedge t_{2}$. Then for any $0\leq t\leq t_{0}$ we have from \eqref{unif:r1}, \eqref{unif:r2}, \eqref{unif:r3},
\eqref{unif:r4}, and \eqref{unif:r5}
\begin{eqnarray*}
&&|r_{H}^{tm}(1,(r,\tilde{0}),(r,\tilde{0}))-r^0_{H}(1,(r,\tilde{0}),(r,\tilde{0}))|\\
&=&|\E_{(r,\tilde{0})}[\tau_{H}^{tm}<1, p^{tm}(1-\tau_{H}^{tm},X^{tm}_{\tau_{H}^{tm}},(r,\tilde{0}))]-
\E_{(r,\tilde{0})}[\tau^0_{H}<1, p^0(1-\tau^0_{H},X^0_{\tau^0_{H}},(r,\tilde{0}))]|\\
&\leq&|\E_{(r,\tilde{0})}[1>\tau_{H}^{tm}>1-\delta_{1},\tau_{H}^{tm}<1, p^{tm}(1-\tau_{H}^{tm},X^{tm}_{\tau_{H}^{tm}},(r,\tilde{0}))]| +\\
&&|\E_{(r,\tilde{0})}[1>\tau^0_{H}>1-\delta_{1},\tau^0_{H}<1, p^0(1-\tau^0_{H},X^0_{\tau^0_{H}},(r,\tilde{0}))]|+\\
&&| \E_{(r,\tilde{0})}\left[T^{tm}\leq1, \tau_{H}^{tm}<1-\delta_{1}, p^{tm}(1-\tau_{H}^{tm},X^{tm}_{\tau_{H}^{tm}},(r,\tilde{0}))\right]| + \\
&&|\E_{(r,\tilde{0})}\left[T^{tm}\leq 1,\tau^0_{H}<1-\delta_{1}, p^0(1-\tau^0_{H},X^0_{\tau^0_{H}},(r,\tilde{0}))\right]| \\
&& +| \E_{(r,\tilde{0})}\left[T^{tm}>1, \tau^0_{H}<1-\delta_{1}, p^{tm}(1-\tau^{0}_{H},X^0_{\tau^0_{H}},(r,\tilde{0}))\right] \\
&&- \E_{(r,\tilde{0})}\left[T^{tm}>1,\tau^0_{H}<1-\delta_{1}, p^0(1-\tau^0_{H},X^0_{\tau^0_{H}},(r,\tilde{0}))\right]| \\
&<&5\eps.
\end{eqnarray*}
\qed

As in \cite{BKS}, we need to divide the Lipschitz open set $D$ into a good set and a bad set.
We recall several geometric facts
about the Lipschitz open set.

\begin{defn}\label{good set}
 Let $\eps, r>0$. We say that $G\subset \partial D$ is $(\eps,r)$-good if for each point $p\in G$, the unit inner normal $\nu(p)$ exists
and
$$
B(p,r)\cap \partial D\subset \{x: |(x-p) \cdot \nu(p)|<\eps |x-p|\}.
$$
\end{defn}
If $G$ is an $(\eps,r)$-good subset of $\partial D$, then using this definition
we can construct a good subset $\mathcal{G}$ of the points near the boundary:
$$
\mathcal{G}=\bigcup_{p\in G}\Gamma_{r}(p,\eps),
$$
where $\Gamma_{r}(p,\eps)=\{x: (x-p)\cdot \nu(p) >\sqrt{1-\eps^{2}}|x-p|\}\cap B(p,r)$.

The next lemma is \cite[Lemma 2.7]{BKS} and it says the measure of the set of the bad points near the boundary is small.
Note that even though \cite[Lemma 2.7]{BKS} is stated for a bounded Lipschitz domain, the proof remains true for
a bounded Lipschitz open set.
\begin{lemma}[Lemma 2.7 in \cite{BKS}]\label{bad set:boundary}
Suppose $\eps \in (0, 1/2)$, $r>0$ and that $G$ is a measurable $(\eps,r)$-good subset of $\partial D$.
There exists $s_{0}(\partial D,G)>0$ such that for all $s<s_{0}$
$$
\left|\{x\in D : \delta_{D}(x)<s\}\setminus \mathcal{G}\right|\leq s\left[\mathcal{H}^{d-1}(\partial D\setminus G)+
\eps\left(3+\mathcal{H}^{d-1}(\partial D)\right)\right].
$$
\end{lemma}

The next lemma is about the existence of a good subset $G\subset \partial D$.
Again the lemma remains true for a bounded Lipschitz open set $D$.

\begin{lemma}[Lemma 2.8 in \cite{BKS}]\label{good set2}
For any $\eps>0$ there exists $r>0$ such that an $(\eps,r)$-good set $G\subset \partial D$ exists and
$$
\mathcal{H}^{d-1}(\partial D\setminus G)<\eps.
$$
\end{lemma}

The two lemmas above imply that
$$
\left|\{x\in D : \delta_{D}(x)<s\}\setminus \mathcal{G}\right|\leq s\eps\left(4+\mathcal{H}^{d-1}(\partial D)\right).
$$

For any $\eps\in (0, 1/4)$, we fix the $(\eps,r)$-good set from Lemma \ref{good set2}
and construct $\mathcal{G}$ from $G$.
We choose $r$ to be smaller than the minimal distances between (finitely many) components of $D$.
For any $x\in \mathcal{G}$, there exists $p(x)\in \partial D$ such that $x\in \Gamma_r(p(x), \eps)$.
Next we define inner and outer cones as follows
\beq\label{inner cone}
I_{r}\left(p(x)\right)=\{y:\left(y-p(x)\right)\cdot \nu\left(p(x)\right)>\eps |y-p(x)|\}\cap B(p(x),r),
\eeq
\beq\label{outer cone}
U_{r}\left(p(x)\right)=\{y:\left(y-p(x)\right)\cdot \nu\left(p(x)\right)<-\eps |y-p(x)|\}\cap B(p(x),r).
\eeq

It follows from \cite[(2.20)]{BKS} that there exists a half-space $H^{*}(x)$ such that
\beq\label{half space}
x\in H^{*}(x), \quad \delta_{H^{*}(x)}(x)=\delta_{D}(x), \quad
I_{r}\left(p(x)\right) \subset H^{*}(x) \subset U_{r}\left(p(x)\right)^{c}.
\eeq

Now we are ready to prove Theorem \ref{RSP:main2}.

\noindent{\bf Proof of Theorem \ref{RSP:main2}.}
Fix $\eps\in (0, 1/4)$, the $(\eps,r)$-good set from Lemma \ref{good set2} and
the $\mathcal{G}$ constructed from $G$.
From the definition of the trace we have
\begin{eqnarray*}
 &&-t^{d/\alpha}\int_{D}r_{D}^{m}(t,x,x)dx=t^{d/\alpha}\int_{D}\left(p_{D}^{m}(t,x,x)-p^{m}(t,x,x)\right)dx\\
&=&t^{d/\alpha}Z_{D}^{m}(t)-t^{d/\alpha}\int_{D}p^{m}(t,x,x)dx\\
&=&t^{d/\alpha}Z_{D}^{m}(t)-t^{d/\alpha}\int_{D}\left(p^{0}(t,x,x)-\left(p^{0}(t,x,x)-p^{m}(t,x,x)\right)\right)dx\\
&=&t^{d/\alpha}Z_{D}^{m}(t)-C_{1}|D|+t^{d/\alpha}\int_{D}\left(p^{0}(t,x,x)-p^{m}(t,x,x)\right)dx.
\end{eqnarray*}
Hence it follows from Lemma \ref{RSP:power terms} that in order to prove Theorem \ref{RSP:main2}
we must show that for given $\eps\in (0, 1/4)$
there exists a $t_{0}>0$ such that for any $0<t<t_{0}$,
$$
\left|t^{d/\alpha}\int_{D}r_{D}^{m}(t,x,x)dx-C_{2}\mathcal{H}^{d-1}(\partial D)t^{1/\alpha}\right|\leq c(\eps)t^{1/\alpha},
$$
where $c(\eps)\rightarrow 0$ as $\eps \rightarrow 0$.
As in the proof of \cite[Theorem 1.1.]{BKS}
we split the region of integration into three sets
\begin{eqnarray*}
&&\mathcal{D}_{1}=\{x\in D\setminus \mathcal{G} : \delta_{D}(x)<s \},\\
&&\mathcal{D}_{2}=\{x\in D\cap \mathcal{G} : \delta_{D}(x)<s \},\\
&&\mathcal{D}_{3}=\{x\in D : \delta_{D}(x)\geq s\},
\end{eqnarray*}
where s must be smaller than the $s_{0}$ given by Lemma \ref{bad set:boundary}. For small enough $t$ we can take
$$
s=t^{1/\alpha}/\sqrt{\eps}\nonumber.
$$
It is shown in \cite[(3.2) and (3.4)]{BKS} that
\beq\label{R1,R3:SSP}
t^{d/\alpha}\int_{\mathcal{D}_{1}\cup\mathcal{D}_{3}}r_{D}^{0}(t,x,x)dx\leq c(\eps)t^{1/\alpha}
\eeq
where $c(\eps)\rightarrow 0$ as $\eps \rightarrow 0$.
Hence it follows from Lemma \ref{r:relation} and \eqref{R1,R3:SSP} that
\beq\label{R1,R3:RSP}
t^{d/\alpha}\int_{\mathcal{D}_{1}\cup\mathcal{D}_{3}}r_{D}^{m}(t,x,x)dx\leq c(\eps)e^{2mt}t^{1/\alpha}.
\eeq

Now we deal with the integral on $\mathcal{D}_{2}$.
Let $H^{*}(x)$, $I_{r}\left(p(x)\right)$, $U_{r}\left(p(x)\right)$ be defined by \eqref{inner cone}, \eqref{outer cone} and \eqref{half space}.
We have
$$
I_{r}\left(p(x)\right)\subset H^{*}(x)\subset U_{r}\left(p(x)\right)^{c}.
$$
Since $r$ is less than the minimal distances between components of $D$, we also have
$$
\quad I_{r}\left(p(x)\right)\subset D\subset U_{r}\left(p(x)\right)^{c}.
$$
Since $I_{r}\left(p(x)\right) \subset U_{r}\left(p(x)\right)^{c}$,
By an argument similar to that used in Lemma \ref{RSP:2nd term 2} we have
\begin{eqnarray}\label{R2:RSP1}
&&\left|r_{D}^{m}(t,x,x)-r_{H^{*}(x)}^{m}(t,x,x)\right|\nonumber\\
&\leq& r_{I_{r}\left(p(x)\right)}^{m}(t,x,x)-r_{U_{r}\left(p(x)\right)^{c}}^{m}(t,x,x)\nonumber\\
&\leq&e^{2mt}\left(r_{I_{r}\left(p(x)\right)}^{0}(t,x,x)-r_{U_{r}\left(p(x)\right)^{c}}^{0}(t,x,x)\right).
\end{eqnarray}
Now it follows from \cite[Proposition 3.1.]{BKS} and \eqref{R2:RSP1} that
\begin{eqnarray*}
&&t^{d/\alpha}\int_{\mathcal{D}_{2}}\left|r_{D}^{m}(t,x,x)-r_{H^{*}(x)}^{m}(t,x,x)\right|dx\\
&\leq& ce^{2mt}\left(\eps^{1-\alpha/2}\vee \sqrt{\eps}\right)\mathcal{H}^{d-1}(\partial D)t^{1/\alpha}
\int_{0}^{\infty}\left(r^{-d-\alpha+1}\wedge 1\right)dr.
\end{eqnarray*}

Finally we will show that the integral
$$
t^{d/\alpha}\int_{\mathcal{D}_{2}}r_{H^{*}(x)}^{m}(t,x,x)dx
$$
gives the second term $C_{2}\mathcal{H}^{d-1}(\partial D)t^{1/\alpha}$ plus an error term of order $c(\eps)t^{1/\alpha}$.
Recall that
$$
r_{H^{*}(x)}^{m}(t,x,x)=f_{H^{*}}^{m}(t,\delta_{H^{*}(x)})=f_{H}^{m}(t,\delta_{D}(x)).
$$
Hence we have
\begin{eqnarray*}
&& t^{d/\alpha}\int_{\mathcal{D}_{2}}r_{H^{*}(x)}^{m}(t,x,x)dx\\
&=&t^{d/\alpha}\int_{\mathcal{D}_{2}}f_{H}^{m}(t,\delta_{D}(x))dx\\
&=&t^{d/\alpha}\int_{D}f_{H}^{m}(t,\delta_{D}(x))dx-t^{d/\alpha}\int_{\mathcal{D}_{1}\cup\mathcal{D}_{3}}f_{H}^{m}(t,\delta_{D}(x))dx.
\end{eqnarray*}
By an argument similar to that used to get \eqref{R1,R3:RSP} we have that
$$
t^{d/\alpha}\int_{\mathcal{D}_{1}\cup\mathcal{D}_{3}}f_{H}^{m}(t,\delta_{D}(x))dx\leq c(\eps)t^{1/\alpha},
$$
where $c(\eps)\rightarrow 0$ as $\eps \rightarrow 0$.
From the (approximate) scaling property of the relativistic stable process, we have
$$
t^{d/\alpha}\int_{D}f_{H}^{m}(t,\delta_{D}(x))dx=\int_{D}f_{H}^{mt}(1,\delta_{D}(x)/t^{1/\alpha})dx.
$$
Now apply Lemmas \ref{conv2} and \ref{conv:unif} to the function $r\rightarrow f_{H}^{mt}(1,r)$ and we get for small enough $t$
$$
\left|\int_{D}f_{H}^{mt}(1,\delta_{D}(x)/t^{1/\alpha})dx-C_{2}\mathcal{H}^{d-1}(\partial D)t^{1/\alpha}\right|\leq\eps t^{1/\alpha}.
$$
\qed

\medskip
{\bf Acknowledgement}
After the first version of this paper, which only contains Theorem \ref{RSP:main}, was finished,
the first named author sent it to Professor Ba\~nuelos. Professor Ba\~nuelos encouraged us
to work out the Lipschitz case. We thank him for his encouragement and for his helpful comments
on a later version of the paper.

\begin{singlespace}

\end{singlespace}
\end{doublespace}

\vskip 0.3truein

{\bf Hyunchul Park}

Department of Mathematics, University of Illinois, Urbana, IL 61801,
USA

E-mail: \texttt{hpark48@illinois.edu}

\bigskip

{\bf Renming Song}

Department of Mathematics, University of Illinois, Urbana, IL 61801,
USA

E-mail: \texttt{rsong@math.uiuc.edu}

\end{document}